\documentclass[11pt]{article}
\usepackage[english]{babel}
\usepackage{amsfonts,amsmath,mathrsfs,stmaryrd,amsthm,alltt,float}
\usepackage{graphicx}
\usepackage{amssymb}
\usepackage{amsmath}
\usepackage{verbatim} 
\usepackage{natbib}
\usepackage{color}
\usepackage{hyperref}
\usepackage[nolists]{endfloat}

\renewcommand{\hat}{\widehat}

\renewcommand{\tilde}{\widetilde}

\def\b#1{\mbox{\boldmath $#1$}}

\def\bSig\mathbf{\Sigma}

\def\b#1{\mbox{\boldmath $#1$}}

\title{\vspace*{-1.5cm} A discrete time event-history approach to informative drop-out in\\
multivariate latent Markov models with covariates\vspace*{5mm}}
\author{Francesco Bartolucci\\
{\small Department of Economics, Finance, and Statistics}\\
{\small University of Perugia (IT)}\\
{\small email: bart@stat.unipg.it}\vspace*{5mm}\and 
Alessio Farcomeni\\
{\small Department of Public Health and Infectious Diseases}\\
{\small Sapienza - University of Rome (IT)}\\
{\small email: alessio.farcomeni@uniroma1.it}}
\date{}
\begin{document}

\label{firstpage}

\maketitle

\begin{abstract}
Latent Markov (LM) models represent an important tool of analysis of longitudinal data when response variables are affected by time-varying unobserved heterogeneity, which is accounted for by a hidden Markov chain. In order to avoid bias when using a model of this type in the presence of informative drop-out, we propose an event-history (EH) extension of the LM approach that may be used with multivariate longitudinal data, in which one or more outcomes of a different nature are observed at each time occasion. The EH component of the resulting model is referred to the interval-censored drop-out, and bias in LM modeling is avoided by correlated random effects, included in the different model components, which follow a common Markov chain. In order to perform maximum likelihood estimation of the proposed model by the Expectation-Maximization algorithm, we extend the usual backward-forward recursions of Baum and Welch. The algorithm has the same complexity of the one adopted in cases of non-informative drop-out. Standard errors for the parameter estimates are derived by using the Oakes' identity. We illustrate the proposed approach through an application based on data coming from a medical study about primary biliary cirrhosis in which there are two outcomes of interest, the first of which is continuous and the second is binary. 
\end{abstract}

\newpage

\section{Introduction}
\label{intro}

In longitudinal studies, subjects may be lost to follow-up due to events, like death, which are associated with the outcome of interest. In these cases, an informative drop-out arises that must be properly modeled in order to avoid biased estimates. From the reverse perspective, the time trend of a longitudinal measurement may predict the risk of an event (e.g., a steadily decreasing CD4 count is predictive of adverse events in HIV patients); see for instance \cite{follmann:95} for a general account of related longitudinal and survival processes. 

A common approach to deal with informative drop-out is via shared-parameter models \citep[e.g.,][]{wu:carrol:88,follmann:95}, where both longitudinal and survival mechanisms are assumed to share a latent Gaussian variable. Semiparametric shared-parameter models have been proposed by \cite{tson:verb:lesa:09}. A discrete random effect, along the lines of this work, is adopted by \cite{roy:03} to deal with an ordinal latent class model. Another approach to modeling informative drop-out is that of \cite{WULF:TSIAT:97} and \cite{riz:10}, where the risk of an event at time $t$ is influenced by the expected value of the longitudinal response at the same time. The resulting Joint Model (JM) includes both the longitudinal fixed and the random effects in the computation of the hazard function. There are very few generalizations of JMs to the case of discrete longitudinal outcomes. Notable exceptions are that of \cite{rizo:ghos:11}, who propose generalized linear models for the longitudinal outcome in a Bayesian framework, and \cite{vivi:alfo:rizo:13}, who deal with a similar model in a classical maximum likelihood framework. 

A limitation of shared-parameter models and JMs is that latent variables, in the form of subject-specific parameters, are time constant. The effect of unobserved heterogeneity, therefore, cannot evolve or must evolve over time along a pre-specified parametric form (e.g., linear through inclusion of a random slope). Further, the use of continuous latent variables may hinder possibilities of clustering subjects with respect to common unobserved heterogeneity and risk of event/drop-out. Finally, when the outcome of interest is categorical, it may be viewed as measuring, with error, a latent {\it discrete} rather than continuous variable. In order to overcome these limitations, latent Markov (LM) models represent a very flexible and convenient way of modeling categorical outcomes which are repeatedly measured over time; see \cite{bart:farc:penn:12} for an overview. The basic assumption of these models is that the response variables, which are longitudinally observed, are conditionally independent given a hidden first-order Markov chain which accounts for the unobserved heterogeneity.

Despite the relevance of LM models, there are very few extensions of these models to deal with drop-out. A model that is somehow in spirit to our approach is that of \cite{albe:00}, which jointly models the outcome and missing mechanism. The latter is assumed to follow a {\it manifest} first-order Markov chain, and the two processes are linked since the outcome is used to model the missingness indicators. In a Bayesian framework, \cite{spag:et:al:11} propose a simple three-state latent chain, in which one state is actually not latent and corresponds to drop-out. Finally, in \cite{maru:13} the time to drop-out is used in a model for the initial and transition probabilities of the hidden Markov chain. 

In this paper, we propose a different approach with respect to the ones mentioned above and that, at least to our knowledge, has not been previously considered in the literature. In the proposed approach, the manifest distribution is jointly referred to the longitudinal and drop-out processes; the corresponding time-varying unobserved heterogeneity structure evolves according to the same initial and transition distributions.  Our approach falls into the class of non-ignorable random-coefficient-based drop-out models as defined in \cite{litt:95}. 

In the proposed approach, the longitudinal outcomes are modeled through generalized linear mixed effects models \citep[e.g.,][]{mccu:sear:04,fitz:garr:lair:ware:04} and a discrete event-history (EH) model \citep[see][and references therein]{stee:11} is used for the drop-out process. It is important to stress that this is an approach to multivariate longitudinal data, in the sense that even more outcomes and of a different nature can be observed at each time occasions. These outcomes are assumed to be conditionally independent given the random effects. A flexible dependence structure is obtained as the random effects are distributed according to a single first-order homogeneous latent Markov chain with a finite number of states. Subjects in the same latent state share class-specific intercepts for the longitudinal models, and a class-specific intercept for the EH model, and can also share common regression coefficients for the covariates. The resulting estimates are easily interpretable, and the model is reasonable as it is natural to expect that longitudinal outcomes and drop-out share the same sources of (time-varying) unobserved heterogeneity, which can have different effects on each of them. 

For the proposed model we perform maximum likelihood estimation by an
Expectation-Maximization (EM) algorithm
\citep{baum:et:al:70,demp:lair:rubi:77}. This requires an extension of
the forward-backward recursions
\citep{baum:et:al:70,welch:2003,bart:farc:penn:12} to account for
informative drop-out. We also pay attention to the computation of the standard errors for the parameter estimates by employing a result due to \cite{oake:99}; see also \cite{bart:farc:13}. Moreover, for
model selection, and in particular for the choice of the number of
latent states, 
we suggest the use of the Bayesian Information Criterion \citep[BIC,][]{sch:78}.

The remainder of the paper is organized as follows. In Section
\ref{sec:model} we illustrate the proposed class of models. In Section
\ref{sec:inference} we illustrate maximum likelihood inference for
these models. Finally, the approach is illustrated in Section
\ref{sec:application} through an application based on primary biliary
cirrhosis data in which there are two outcomes of interest for each time occasions, one of which is continuous and the other is binary. We provide some concluding remarks in Section
\ref{concl}. 

\section{A class of latent Markov models with informative drop-out}\label{sec:model}

We consider a longitudinal study on a sample of $n$ subjects, or more generally sample units, in which $s$ follow-up time occasions are scheduled for each subject. We assume that $s$ is known in advance and equal for all subjects. A generalization to a subject-specific number of occasions is straightforward, while if $s$ is not known in advance all proposed inference can be thought of as being conditional on the maximum number of follow-up times which is observed. 

For every subject $i$ and time occasion $t$, with $i=1,\ldots,n$ and $t=1,\ldots,T$, we observe $r$ response
variables, denoted by $Y_{hit}$, $h=1,\ldots,r$; we also denote by
$T_i$ the last time occasion of observation for subject $i$, so that
drop-out occurs before occasion $T_i+1$. As will be clear in the
following, also $T_i$ is a random variable, the distribution of which
depends on observable covariates and latent variables for the
unobserved heterogeneity. A realization of the $h$-th response variable is
denoted by $y_{hit}$ and, accordingly, a realization of $T_i$ is denoted by $t_i$. The observed outcomes for the same subject $i$ and time occasion $t$ are collected in the column vector $\b y_{it}=(y_{1it},\ldots,y_{rit})'$. Also note that if there is drop-out then $t_i<s$, whereas $t_i=s$ indicates that a complete record of outcomes is observed for the $i$-th subject.

Let $D_{it}$, $i=1,\ldots,n$, $t=1,\ldots,s$, denote a binary
random variable equal 1 if subject $i$ drops out from the study
after occasion $t$ and before occasion $t+1$, that is, $T_i=t$, and to 0
otherwise. The basic assumption of the proposed model is that, given a
discrete latent variable $U_{it}$ with $k$ support points and the
available covariates, the response variables $Y_{1it},\ldots,Y_{rit}$
are conditionally independent and they are also independent of
$D_{it}$. In particular, we denote by $\b x_{hit}$ the column vector
of covariates affecting $Y_{hit}$ and by $\mu_{hit}(u)$ the expected
value of this response variable given these covariates and
$U_{it}=u$. Note that the dependence on the covariates of this mean is
not explicitly indicated since these covariates are considered as
fixed and given; this convention will be used for the notation throughout the
paper. Similarly, we denote by $\b z_{it}$ the column vector of
covariates affecting $D_{it}$ and by $p_{it}(u)$ the probability that
$D_{it}=1$ given these covariates and $U_{it}=u$. Then, for
$u=1,\ldots,k$, we assume 
\begin{eqnarray*}
\begin{cases}
g_1[\mu_{1it}(u)] = \alpha_{1u} + \b x_{1it}'\b\beta_1,\qquad \mbox{$t=1,\ldots,s$,} \\
\quad\quad\vdots\quad\quad\quad\quad\vdots\\
g_r[\mu_{rit}(u)] = \alpha_{ru} +\b x_{rit}' \b\beta_r, \qquad \mbox{$t=1,\ldots,s$,} \\
g[p_{it}(u)] =  \gamma_u + \b z_{it}'\b\delta,\qquad\qquad\ 
\mbox{$t=1,\ldots,s-1$,}\\
p_{is}(u) = 1, 
\end{cases}
\end{eqnarray*}
where $g_h(\cdot)$ and $g(\cdot)$ are appropriate link functions and every $Y_{hit}$ is assumed to have conditional distribution belonging to the regular exponential family \citep{mcc:nel:89}. The reason why $p_{is}(u) = 1$ is that the follow-up surely stops after occasion $s$ for all subjects. 

In the dataset used to illustrate the proposed approach (see Section \ref{sec:application}) there are $r=2$ response variables. The first of these variables is continuous and the second is binary. Hence, for this application we choose $g_1(\cdot)$ as corresponding to the identity link and $g_2(\cdot)$ as corresponding to the logit link; for the distribution of the first variable we also have a dispersion parameter indicated, in general, by $\sigma_h^2$. The model for each longitudinal outcome is a classical generalized linear mixed effects model, while the time to drop-out follows a geometric distribution as in classical discrete time EH models \citep[e.g.,][]{stee:11}. We regard drop-out as a trial within each time interval. It is shown in \cite{stee:diam:amin:96} that the resulting likelihood is that of a Bernoulli; consequently, we have
$$
\Pr(T_i=t_i|U_{i1}=u_1,\ldots,U_{it_i}=u_{t_i},\b z_{i1},\ldots,\b z_{it_i}) = 
p_{it_i}(u_{t_i})\prod_{t=1}^{t_i-1}[1-p_{it}(u_t)] 
$$
where, since $p_{is}(u)=1$, the probability $p_{it_i}(u_i)$
disappears from the above expression when there is no drop-out until
the end of the study ($t_i=s$). This recovers the {\it truncated}
geometric distribution given that $s$ is finite.
In order to derive the results in Section \ref{sec:inference} it is also important to note that
\[
\Pr(T_i>t|U_{it}=u) = \prod_{j=1}^t[1-p_{ij}(u)],\quad t=1,\ldots,s,
\]
and that
\[
\Pr(T_i=t_i|T_i>t,U_{i,t+1}=u_{t+1},\ldots,U_{it_i}=u_{t_i})=p_{it_i}(u_{t_i})\prod_{j=t+1}^{t_i-1}[1-p_{ij}(u_j)],\quad t=1,\ldots,t_i-1.
\]

Note that the only information that is used in the model is the interval censored event time (i.e., that drop-out occurs between $t_{i}$ and $t_{i+1}$). We further note that when continuous-time durations are grouped into discrete intervals, a continuous-time hazard model (with constant hazard within each interval) would lead to the complementary log-log model with $g[p_{it}(u)] = \log\{-\log[1-p_{it}(u)]\}$; see \cite{kalb:pren:02}. In our application we use a logit link, which we find more convenient.  

Regarding the distribution of the latent variables, we assume that, for $i=1,\ldots,n$, the sequence $U_{i1},\ldots,U_{is}$ follows a Markov chain with initial probabilities 
$$
\pi_u=\Pr(U_{i1}=u),\quad u=1,\ldots,k, 
$$
which are collected in the column vector $\b\pi$, and time-homogeneous transition probabilities
$$
\pi_{uv}=\Pr(U_{it}=v|U_{i,t-1}=u),\quad u,v=1,\ldots,k, 
$$
which are collected in the transition matrix $\b\Pi$. The parameters of the model that must be estimated are then $\b\pi$, $\b\Pi$, $\b\alpha$, $\b\beta$, $\b\gamma$, $\b\delta$, where $\b\alpha$ and $\b\gamma$ are column vectors with elements $\alpha_{hu}$, $h=1,\ldots,r$, $u=1,\ldots,k$, and $\gamma_u$, $u=1,\ldots,k$, respectively. Similarly, $\b\beta$ is the vector obtained by casting the vectors $\b\beta_h$, $h=1,\ldots,r$. 

The degree of dependence between the longitudinal and the survival
processes is measured by the total variation of the support points of
the EH model. Consequently, when $k=1$ then drop-out is
non-informative. A classical LM model with $k>1$ states, where
drop-out is non-informative, may obviously also be obtained by fixing
$\gamma_u=\gamma$, $u=1,\ldots,k$. There is sensitivity to drop-out as
soon as $\gamma$ is not constant, 
hence we can use a formal test based on a likelihood ratio statistic 
for the null hypothesis $H_0:\gamma_1=\cdots=\gamma_k$.

Note that the link between the random effects $\gamma_u$ and $\alpha_{hu}$ is based on the assumption that they follow the same latent Markov rule. This approach is similar in spirit to situations in which a copula is
used to model the dependence of random effects of the two processes,
like in \cite{rizopoulos:08b}.

\section{Likelihood inference}\label{sec:inference}
We begin considering the observed likelihood 
$$
L(\b\theta) = \prod_{i=1}^n
f(\b y_{i1},\ldots,\b y_{it_i},t_i),
$$
where $\b\theta$ is a short-hand notation for all the model parameters and $f(\b y_{i1},\ldots,\b y_{it_i},t_i)$ is the density or probability of the observed outcomes given the covariates $\b x_{it}$ and $\b z_{it}$ until occasion $t_i$; these covariates are not explicitly indicated as they are fixed and given. This component of the observed likelihood can be expressed as
\begin{equation}
\label{indobs}
\hspace{-0.5cm}f(\b y_{i1},\ldots,\b y_{it_i},t_i)=\sum_{u_1=1}^k\cdots\sum_{u_{t_i}=1}^k
\left(\pi_{u_1}\prod_{t=2}^{t_i} \pi_{u_{t-1}u_t}\right)
\left[\prod_{h=1}^r\prod_{t=1}^{t_i} f(y_{hit}|u_t)\right]
\left\{p_{it_i}(u_{t_i})\prod_{t=1}^{t_i-1} [1-p_{it}(u_t)]\right\},
\end{equation} 
where $f(y_{hit}|u)$ refers to the conditional density or probability of $Y_{hit}$ evaluated at $y_{hit}$, given $U_t=u$ the corresponding covariates, and
\[
p_{it}(u)=\frac{\exp(\gamma_u+\b z_{it}'\b\delta)}{1+\exp(\gamma_u+\b z_{it}'\b\delta)},
\]
for $t=1,\ldots,s-1$ with $p_{is}(u)=1$. 

Expression (\ref{indobs}) can be efficiently computed by an extension of the {\em forward recursion}, which is well known in the hidden Markov literature \citep{baum:et:al:70,zucc:macd:09,bart:farc:penn:12}. First of all, ruling out the trivial case in which $s=1$, we consider the following density or probability for $t=1,\ldots,t_i$ and $u=1,\ldots,k$:
$$
a_{it}(u) =  \left\{\begin{array}{ll} 
f(\b y_{i1},\ldots,\b y_{it},T_i>t,U_{it}=u),&\mbox{if }\:t<t_i,\\
f(\b y_{i1},\ldots,\b y_{it},T_i=t,U_{it}=u),&\mbox{if }\:t=t_i. 
\end{array}\right.
$$
Then, for $t=1$ we have that 
$$
a_{i1}(u) = \left\{\begin{array}{ll} 
\pi_u\left[\prod_{h=1}^rf(y_{hi1}|u)\right][1-p_{it}(u)] & \mbox{if }
t_i>1,\\
\pi_u\left[\prod_{h=1}^rf(y_{hi1}|u)\right]p_{it}(u) & \mbox{if }
t_i=1,
\end{array}\right.
$$
whereas, provided that $t_i>1$, for $t=2,\ldots,t_i$ we have 
$$
a_{it}(v) = \left\{\begin{array}{ll}
\sum_{u=1}^k a_{i,t-1}(u)\pi_{uv}\left[\prod_{h=1}^r f(y_{hit}|v)\right]
[1-p_{it}(v)], & \mbox{if } t<t_i,\\
\sum_{u=1}^k a_{i,t-1}(u)\pi_{uv}\left[\prod_{h=1}^r f(y_{hit}|v)\right]
p_{it}(v), & \mbox{if } t=t_i,
\end{array}\right.
$$
where $v=1,\ldots,k$. At the end of the recursion ($t=t_i$), we have that
$$
f(\b y_{i1},\ldots,\b y_{it_i},t_i)=\sum_{u=1}^k a_{it_i}(u).
$$

In order to maximize the likelihood $L(\b\theta)$, we implement a version of the EM algorithm \citep{demp:lair:rubi:77}, which is based on the {\em complete data likelihood}. Let $w_{it}(u)$ denote a dummy variable equal to 1 if the $i$-th subject is in latent state $u$ at the $t$-th occasion and let $z_{it}(u,v)=w_{i,t-1}(u)w_{it}(v)$ be a dummy variable equal to 1 if there is a transition from latent state $u$ to latent state $v$ at occasion $t$. The logarithm of this likelihood has the following expression: 
\begin{eqnarray}
\ell_c(\b\theta)&=&\sum_{i=1}^n\Bigg(
\sum_{u=1}^k w_{i1}(u) \log\pi_u + \sum_{t=2}^{t_i}\sum_{u=1}^k\sum_{v=1}^k
z_{it}(u,v) \log\pi_{uv}+\sum_{h=1}^r\sum_{t=1}^{t_i}\sum_{u=1}^k w_{it}(u) \log
f(y_{hit}|u)\nonumber\\
&+&\sum_{u=1}^k\left\{\sum_{t=1}^{t_i-1} w_{it}(u)
\log[1-p_{it}(u)]+w_{it_i}(u)
\log p_{it_i}(u)\right\} \Bigg),\label{lc}
\end{eqnarray}
where the second summand involving $z_{it}(u,v)$ disappears if $t_i=1$. 

With reference to our application based on $r=2$ response variables,
the first of which has normal conditional distribution parametrized by
an identity link function and the second of which has Bernoulli
distribution parametrized by a logit link function, we have that
\begin{eqnarray*}
\sum_{h=1}^r\sum_{t=1}^{t_i}\sum_{u=1}^k w_{it}(u) \log
f(y_{hit}|u)&=&\sum_{t=1}^{t_i}\sum_{u=1}^k w_{it}(u)
\log\left[\frac{1}{\sigma_1}\phi\left(\frac{y_{1it}-\alpha_{1u}-\b
    x_{1it}'\b\beta_1}{\sigma_1}\right)\right]\\
&&+\sum_{t=1}^{t_i}\sum_{u=1}^k w_{it}(u) \log
\frac{\exp[y_{2it}(\alpha_{2u}+\b x_{2it}'\b\beta_2)]}
{1+\exp(\alpha_{2u}+\b x_{2it}'\b\beta_2)},
\end{eqnarray*}
where $\phi(\cdot)$ denotes the standard normal density function.

The EM algorithm alternates two steps until convergence: first, the
conditional expected value of the complete data log-likelihood is
obtained (E step). The resulting expression is then maximized with
respect to $\b\theta$ (M step). 

The EM is guaranteed to converge to a local optimum of the observed likelihood. In order to increase the chances of reaching the global maximum, we use a multistart strategy. The first initialization is based on estimating $\b\beta$ and $\b\delta$ using separate generalized linear models. The parameters $\b\alpha$ and $\b\gamma$ are then centered on the intercepts obtained by the previous generalized linear models. Finally, $\b\Pi$ is initialized so that all off-diagonal elements $\pi_{uv}$, with $v\neq u$, are equal to $1/[k(k-1)]$. The other starting solutions are obtained by randomly perturbing the parameter estimates obtained at convergence from the deterministic one. 

Finally, we use the value of 
the likelihood at convergence to compute the index on which the
Bayesian Information Criterion (BIC) for model choice is based \citep{sch:78}. 

\subsection{E step}

At the E-step, the conditional expected value of (\ref{lc}) is simply
computed  by a plug-in of the expected values of $w_{it}(u)$ and
$z_{it}(u,v)$ given the observed data and the current value of the
parameters. These expected values are denoted by $\tilde{w}_{it}(u)$
and $\tilde{z}_{iu}(u,v)$, and can be computed by means of an appropriate backward recursion adapted from the hidden Markov literature as we illustrate below.
First of all, we consider the probabilities:
\[
b_{it}(u) = f(\b y_{i,t+1},\ldots,\b y_{it_i},T_i=t_i|T_i>t,U_{it}=t),\quad t=1,\ldots,t_i.
\]
We have that $b_{it_i}(u)=1$ for $u=1,\ldots,k$, whereas, for $t=1,\ldots,t_i-1$ and $u=1,\ldots,k$, we have
\[
b_{it}(u)=
\left\{ \begin{array}{ll} 
\sum_{v=1}^k b_{i,t+1}(v)\pi_{uv} \left[\prod_{h=1}^rf(y_{hi,t+1}|v)\right]
[1-p_{i,t+1}(u)],& t<t_i-1,\\
\sum_{v=1}^k b_{i,t+1}(v)\pi_{uv} \left[\prod_{h=1}^rf(y_{hi,t+1}|v)\right]
p_{i,t+1}(u),&  t=t_i-1,
\end{array}\right.
\]
where as before $p_{is}(u)=1$. 
At the end of the recursion, we obtain the following expected values:
$$
\tilde{w}_{it}(u) = \frac{a_{it}(u)b_{it}(u)}{f(\b y_{i1},\ldots,\b y_{it_i},t_i)},
\quad t=1,\ldots,t_i,\:u=1,\ldots,k;
$$
provided that $t_i\geq 2$, we also obtain, for $u,v=1,\ldots,k$, the expected values
$$
\tilde{z}_{it}(u,v) = \frac{a_{i,t-1}(u)\pi_{uv}\left[\prod_{h=1}^rf(y_{hit}|v)\right][1-p_{it}(v)]b_{it}(v)}{f(\b y_{i1},\ldots,\b y_{it_i},t_i)},\quad t=2,\ldots,t_i-1,
$$
with
$$
\tilde{z}_{it_i}(u,v) = \frac{a_{i,t_i-1}(u)\pi_{uv}\left[\prod_{h=1}^rf(y_{hit_i}|v)\right]p_{it_i}(v)b_{it_i}(v)}{f(\b y_{i1},\ldots,\b y_{it_i},t_i)}.
$$

\subsection{M step}

At the M-step the conditional expected value of (\ref{lc}) is
maximized by separately maximizing its addends. It is straightforward
to check that explicit solutions are available for $\pi_u$; in particular we have
\[
\pi_u=\frac{1}{n}\sum_{i=1}^n \tilde{w}_{i1}(u),\quad u=1,\ldots,k.
\]
In a similar way, for the transition probabilities we have the following solution: 
\[
\pi_{uv}=\frac{\sum_{i=1}^n\sum_{t=2}^{t_i} \tilde{z}_{it}(u,v)}
{\sum_{i=1}^n\sum_{t=2}^{t_i} \tilde{w}_{i,t-1}(u)},\quad u,v=1,\ldots,k, 
\]
where the outer sum at numerator and denominator is extended to all $i$ such that $t_i\geq2$.

For what concerns the other parameters, we can use $r+1$ separate Newton-Raphson algorithms, similar to that used for standard generalized linear models. The first $r$ algorithms are used to update the parameters $\alpha_{hu}$ and $\b\beta_h$ by maximizing the following expression: 
$$
\sum_{i=1}^n\sum_{t=1}^{t_i}\sum_{u=1}^k \tilde{w}_{it}(u) \log
f(y_{hit}|u), 
$$
for $h=1,\ldots,r$. For instance, in the case of a binary outcome with logit link, this corresponds to the maximization of 
$$
\sum_{i=1}^n\sum_{t=1}^{t_i}\sum_{u=1}^k \tilde{w}_{it}(u) \log
\frac{\exp[y_{hit}( \alpha_{hu}+\b x_{hit}'\b\beta_h)]}
{1+\exp(\alpha_{hu}+\b x_{hit}'\b\beta_h)},
$$
whereas explicit solutions are available with Gaussian outcomes also for the dispersion parameter $\sigma_h^2$. The last Newton-Raphson is used to maximize, with respect to $\b\gamma$ and $\b\delta$, the expression 
$$
\sum_{i=1}^n\sum_{u=1}^k\left\{\sum_{t=1}^{t_i-1} \tilde{w}_{it}(u)\log[1-p_{it}(u)]+
\tilde{w}_{it_i}(u)\log p_{it_i}(u)\right\}.
$$
\subsection{Computation of standard errors}
In order to estimate the standard errors for the parameter estimates, we use an approach based on
\cite{oake:99}'s equality, simplified by the fact that the score vector for the 
expected complete log-likelihood is often available in a closed form,
given that it corresponds to the score of $r+1$ generalized linear models; see also \cite{bart:farc:13}. 
For instance, for binary outcomes with 
logit link, the first derivative of the expected complete log-likelihood 
with respect to $\xi_{hu}$ is easily derived as: 
$$
\sum_{i=1}^n \sum_{t=1}^{t_i} \tilde{w}_{it}(u) [y_{hit}-
\mu_{hit}(u)],
$$
whereas the first derivative with respect to $\b\beta_h$ can be derived as 
$$
\sum_{i=1}^n \sum_{t=1}^{t_i}\sum_{u=1}^k  \tilde{w}_{it}(u)
[y_{hit}-\mu_{hit}(u)]\b x_{hit}. 
$$
Moreover, the first derivative of the expected complete log-likelihood 
with respect to $\gamma_u$ can be derived as: 
$$
\sum_{i=1}^n\left\{-\sum_{t=1}^{t_i-1} \tilde{w}_{it}(u)p_{it}(u)+
\tilde{w}_{it_i}(u)[1-p_{it_i}(u)]\right\},
$$
and finally the derivative with respect to $\b\delta$ is 
$$
\sum_{i=1}^n\sum_{u=1}^k\left\{-\sum_{t=1}^{t_i-1} \tilde{w}_{it}(u)p_{it}(u)+
\tilde{w}_{it_i}(u)[1-p_{it_i}(u)]\right\}\b z_{it}. 
$$
These derivatives are the same that are used in the Newton-Raphson algorithm to implement the M-step as described above.

Once we have the score vector, it can be noted as in \cite{oake:99}
that the observed information is equal to the Jacobian of the
score vector with respect to the parameters, plus the Jacobian of the
score vector with respect to $\tilde{w}_{it}(u)$, seen as a function
of the parameters. For convenience, the second Jacobian 
is computed numerically. 

\section{Application to primary biliary cirrhosis data}
\label{sec:application}

We illustrate the proposed approach by an application based on data coming from a randomized study for treatment of primary biliary cirrhosis. These data were previously analyzed by \cite{murt:et:al:94} and \cite{rizo:verb:mole:10} from a slightly different perspective than the present one. In this study, $n=312$ patients were randomized to a placebo or a treatment based on D-penicillamine. We are interested in evaluating the effect of treatment on a continuous outcome ({\em serum Bilirubin in mg/dl}, $Y_1$) and a binary outcome ({\em presence of edema}, $Y_2$), after adjusting for
certain covariates ({\em drug}, {\em age}, {\em gender}, {\em albumin in gm/dl}, {\em alkaline phosphatase in U/liter}, and {\em SGOT in U/ml at the first visit}) and drop-out. Continuous covariates are zero centered, as the sample mean has been subtracted from each of these covariates. In this application it is reasonable to expect drop-out to be informative as it may be due to death related to the illness or to a transplant. 

The maximum number of follow-up time occasions is $s=16$ for these data and only 1\% of patients have a complete record; the median time to drop-out is $5$. Note that a high serum Bilirubin is very likely to speed the cirrhosis up, hence making drop-out due to death more likely; this enforces our idea that the drop-out cannot be considered as non-informative. In this regard, Table \ref{obsYD} reports the proportion of subjects having a certain number of observations, that is, $\sum_{i=1}^n I(T_i=t)/n$, $t=1,\ldots,T$, the corresponding \cite{kaplan:meier:1958}'s estimates, and the mean of the two outcomes based on the number of survivors.

\begin{table}[ht!]
\caption{\em Observed mean of serum Bilirubin and 
proportion of subjects with edema by time for the primary biliary cirrhosis data, together with the proportion of survivors ($\sum_{i=1}^n I(T_i=t)/n$) and corresponding \cite{kaplan:meier:1958}'s estimates (KM$_t$).}
\label{obsYD}
\begin{center}
\begin{tabular}{cccccc}
\hline\hline
time && $\sum_{i=1}^n I(T_i=t)/n$ & KM$_t$ & mean($Y_1$) & mean($Y_2$) \\
\hline 
1 &&  0.09 & 0.09 & 3.22 & 0.21\\ 
  2 && 0.08 & 0.09 & 3.07 & 0.21 \\ 
  3 && 0.10 & 0.10 & 3.45 & 0.24  \\ 
  4 && 0.14 & 0.14 & 4.26 & 0.27\\ 
  5 && 0.10 & 0.10 & 3.62 & 0.28 \\ 
  6 && 0.07 & 0.07 & 3.91 & 0.30 \\ 
  7 && 0.06 & 0.06 & 3.77 & 0.37 \\ 
  8 && 0.07 & 0.07 & 3.93 & 0.36 \\ 
  9 && 0.05 & 0.05 & 4.01 & 0.37 \\ 
  10 && 0.07 & 0.07 & 3.49 & 0.41\\ 
  11 && 0.05 & 0.05 & 5.14 & 0.40 \\ 
  12 && 0.04 & 0.04 & 4.23 & 0.34 \\ 
  13 && 0.02 & 0.02 & 5.06 & 0.40 \\ 
  14 && 0.02 & 0.02 & 4.32 & 0.43\\ 
  15 && 0.02 & 0.02 & 6.28 & 0.33\\ 
  16 && 0.01 & 0.02 & 5.17 & 0.67\\ 
\hline 
\end{tabular}
\end{center}
\end{table}

We begin the analysis of these data by selecting the number of latent states of the hidden Markov chain. In Table \ref{choosek} we report the maximum log-likelihood, number of parameters, and BIC \citep{sch:78} for increasing values of $k$ from 1 to 4. 

\begin{table}[ht!]
\caption{\em Maximum log-likelihood, number of parameters and BIC
  for different values of $k$ for the primary biliary cirrhosis data.}
\label{choosek}
\begin{center}
\begin{tabular}{cccc}
\hline\hline
$k$ & log-lik. & \# par & BIC \\
\hline
1 & -8290.7 & 28 & 16742.2 \\
2 & -7265.3 & 34 & 14725.9\\
3 & -7254.9 & 42 & 14751.1 \\
4 & -7236.0 & 52 & 14770.6 \\\hline 
\end{tabular}
\end{center}
\end{table}

On the basis of the results in Table \ref{choosek}, we select $k=2$, as the corresponding model has the lowest BIC. In fact, for this model we have maximum log-likelihood of -7265.3 with 34 parameters and then BIC=14725.9. Sensitivity to drop-out can be evaluated by fitting the same model with $k=2$, but with constant latent intercepts (i.e., $H_0:\gamma_1=\gamma_2$). The latter model has a maximum log-likelihood of -7334.5 with 33 parameters and then BIC=14858.5, which is much higher than the previous value. Then, we do not consider hypothesis $H_0$ to be plausible at least when $k=2$, as the likelihood ratio test statistic for this hypothesis, which is equal to 138.4, confirms. The conclusion is that there is evidence of informative drop-out with these data. 

In Table \ref{pars} we report the parameter estimates obtained with our model, and the corresponding LM 
model that assumes the drop-out to be non-informative (including $H_0$). In this table, the intercepts are obtained by averaging the support points. 

\begin{table}[ht]
\caption{\em Parameter estimates and Wald test statistics for the primary biliary cirrhosis data obtained with the proposed LM model and with the corresponding model assuming non-informative drop-out ($H_0:\gamma_1=\gamma_2$). The number of latent states is fixed at $k=2$. The intercept parameters are obtained by averaging the random intercepts.}
\label{pars}
\begin{center}
\begin{tabular}{crrcrr}
\hline\hline
& \multicolumn2c{Proposed Model} && \multicolumn2c{Under $H_0$} \\
\hline 
& \multicolumn5c{\underline{\em Slopes for the serum Bilirubin}}\\
Parameter & Estimate & $t$-statistic &&  Estimate & $t$-statistic \\
\hline 
Intercept & 5.63 & -  && 5.71 & - \\
Treatment & 0.06 & 0.07 && 0.02 & 0.04 \\
Age &  -0.06 & -5.71 && -0.06 & -4.1 \\
Gender (F) & -1.03 & -3.30 && -1.00 &  -3.17 \\
albumin &  -1.31 & -5.61 && -1.47  & -4.84 \\
alkaline ph.$/100$ & 0.03 & 4.74 && 0.02 & 4.25 \\ 
SGOT & 0.02 & 21.01 && 0.02 & 16.92 \\
time & 0.11 &  1.30 && 0.14 & 1.63 \\  
time$^2$ &  -0.00 & -0.12 && -0.00 & -0.33\\
\hline
& \multicolumn5c{\underline{\em Log-odds ratios for the probability of edema}}\\
Parameter & Estimate & $t$-statistic &&  Estimate & $t$-statistic \\
\hline 
Intercept & -0.83 & -  && -0.60 & - \\
Treatment & -0.57 & -3.04 && -0.67 & -2.97 \\
Age &  0.08 & 21.88 && 0.09 & 17.31 \\
Gender (F) &  1.79 & 14.15 && 1.79 &  14.29 \\
albumin &  -2.62 & -22.66 && -2.91  & -21.08 \\
alkaline ph.$/100$ & 0.02 & 10.18 && 0.02 &  9.27 \\ 
SGOT & -0.02 & -24.71 && -0.02 & -23.21 \\
time & 0.58 &  14.24 && 0.67 & 16.01  \\  
time$^2$ &  -0.02 & -6.73 &&  -0.03 & -7.84 \\
\hline
& \multicolumn5c{\underline{\em Log-odds ratios for the probability of drop-out}}\\
Parameter & Estimate & $t$-statistic &&  Estimate & $t$-statistic \\
\hline 
Intercept & -1.06 & -  && -1.62 & - \\
Treatment & 0.01 & 0.24 && 0.09 & 7.41 \\
Age &  0.01 & 0.95 && 0.02 & 2.77 \\
Gender (F) &  0.18 & 9.93 && -0.13 & -3.88 \\
albumin &  -0.91 & -14.83 &&  0.23  & 8.31 \\
alkaline ph.$/100$ & -0.01 & -11.24 && -0.00 & -0.68 \\ 
SGOT & 0.01 & 4.16 && 0.00 & 1.09 \\
time & 0.13 &  6.90 && -0.01 & -8.89 \\  
time$^2$ &  0.01 & 1.09 && 0.00 & 4.13 \\
\hline 
\end{tabular}
\end{center}
\end{table}

From the left panel of Table \ref{pars} (results under the proposed model) we can draw the following conclusions. First of all, the treatment does not seem to be effective on time to drop-out or on level of serum Bilirubin, but at least it significantly decreases the probability of edema, with an OR of 0.57 after adjusting for other factors. Note that, instead, \cite{rizo:verb:mole:10} concluded that the treatment is effective on serum Bilirubin. In our interpretation, this difference is due to the different way of modeling unobserved heterogeneity, that in the approach of \cite{rizo:verb:mole:10} is assumed to increase linearly in time. As far as the other predictors are concerned, it can be seen that with age the level of serum Bilirubin decreases and the probability of edema increases, that females are at higher risk of drop-out and edema but have lower levels of Bilirubin, and that albumin is protective.  

A comparison between the left panel and the right panel of Table \ref{pars} allows us to study the sensitivity to drop-out. Assuming that the drop-out is non-informative (under $H_0$) leads to certain differences in the parameter
estimates and standard errors, especially for the time-event history part for which we may observe that several parameter estimates change sign under $H_0$ with respect to the initial model specification. 

We finally consider the estimated latent distribution parameters. In Table \ref{latent} we report the differences between the latent intercepts and their averages (which are reported in Table \ref{pars}), the initial parameter vector, and the transition matrix. 

\begin{table}[ht]
\caption{\em Differences between the latent intercepts and their average, 
initial parameter vector and hidden transition matrix for the primary
biliary cirrhosis data obtained with the proposed LM model. The number of latent states is fixed at $k=2$.}
\label{latent}
\begin{center}
\begin{tabular}{crrr|c|cc}
\hline\hline\vspace*{-0.32cm}&&&&&&\\
$u$ & \multicolumn1c{$\hat{\alpha}_{1u}$} & \multicolumn1{c}{$\hat{\alpha}_{2u}$} & \multicolumn1{c|}{$\hat{\gamma}_u$} & $\hat{\pi}_u$ & $\hat{\pi}_{u1}$ & $\hat{\pi}_{u2}$  \\\hline
1 & 3.33 & -4.63   & -1.55 & 0.76 & 0.92 & 0.08 \\
2 & 7.94 &  2.97 & -0.57 & 0.24 & 0.06 & 0.94 \\
\hline
\end{tabular}
\end{center}
\end{table}

We observe that the two groups are very well separated. Three out of
four patients start in group 1, which has a slightly lower serum
Bilirubin, propensity to edema and probability of drop-out. This group
is highly persistent, but still it can be seen from the transition
matrix that the proportion of subjects in this group steadily
decreases over time. At $t=5$, only 60\% are in this group, and at
$t=15$ this proportion has decreased to 45\%. As a consequence, it can
be said that the global health status of patients has worsened over
time more than as predicted by the trend estimated by fixed effects (Table \ref{pars}). In particular, the second group is made instead of patients with a sensibly higher Bilirubin, much higher
propensity of edema and slightly higher risk of drop-out. 

\section{Conclusions}
\label{concl}

In this paper, we propose an event-history approach \citep{stee:11}
to account for informative drop-out in LM models for
multivariate longitudinal outcomes. The main features of these models
are that they assume conditional independence among response variables
given a hidden Markov chain. In order to deal with the models formulated in this way, we
extend the usual forward-backward recursions
\citep{baum:et:al:70,demp:lair:rubi:77} to informative drop-out and
consequently propose an extended version of the Expectation-Maximization algorithm
\citep{demp:lair:rubi:77} for maximum likelihood estimation of LM models. 

One of the main features of an LM model with covariates is that it
allows for a form of unobserved heterogeneity that 
is not restricted to be time constant. This form of heterogeneity is
typically explained on the basis of missing covariates, 
and its evolution is accounted for by the inclusion
of random intercepts in the model.
These intercepts are time and individual specific and follow a Markov
chain. Assuming time-constant unobserved heterogeneity may be restrictive, especially when individuals are followed up for many occasions and/or for a long period of time. The advantages of the proposed type of modeling are shown by an application based on a study about primary biliary cirrhosis in which there are two outcomes of interest, the first of which is continuous and the second is binary.

An important assumption of the adopted LM model is that the response
variables referred to the same time occasion are conditionally
independent given the corresponding latent variable and covariates;
this assumption is usually known as {\em local independence}. 
This assumption simplifies model estimation substantially and 
may be explained by considering that all factors affecting
the responses, that is, observed and unobserved covariates, are
properly accounted for. In any case, this form of local independence
can be relaxed by using a multivariate link function, as the one
adopted in the LM framework proposed by \cite{bart:farc:09}.

\section*{Acknowledgments}
Francesco Bartolucci acknowledges the financial support from the grant 
RBFR12SHVV of the Italian Government (FIRB project ``Mixture and latent
variable models for causal inference and analysis of socio-economic data").

\bibliography{biblio}
\bibliographystyle{metron}

\label{lastpage}
\end{document}